\documentclass[11pt]{article}
\usepackage{amsmath}
\usepackage{amsthm,amsfonts,amssymb,epsfig,graphics}
\usepackage[latin1]{inputenc}
\usepackage[english,french]{babel}
\usepackage{hyperref}
\usepackage{pstricks,pst-plot}

\def\homeo{hom\'eomorphisme}

\def\inte{\mathrm{Int}}

\newcommand{\bbR}{{\mathbb{R}}}

\newcommand{\bbD}{{\mathbb{D}}}

\newcommand{\bbS}{{\mathbb{S}}}
\def\bbT{{\mathbb{T}}}

    \def\cV{{\cal V}}
    
\def\cF{{\cal F}}

\def\homeo{\mathrm{Homeo}}
\def\diffeo{\mathrm{Diffeo}}
\def\area{\mathrm{Area}}
\def\diff{{\scriptscriptstyle \mathrm{diff}}}

\def\?{$^{***}$\marginpar{?}}

\newtheorem{theo}{Theorem}
\newtheorem{ques}{Question}
\newtheorem*{ques*}{Question}
\newtheorem*{prop*}{Proposition}
\newtheorem*{conj*}{Conjecture}
\newtheorem*{theo*}{Th\'eor\`eme}
\newtheorem{coro}{Corollary}[section]

\newtheorem{affi*}{Affirmation}
\newtheorem{prop}[coro]{Proposition}
\newtheorem{lemm}[coro]{Lemma}
\def\?{\footnote{?}}

\newlength{\espaceavantspecialthm}
\newlength{\espaceapresspecialthm}
{\vskip \espaceapresspecialthm}

{\vskip \espaceapresspecialthm}

\author{Fr\'ed\'eric Le Roux\footnote{This work was partially supported by the ANR Grant ``Symplexe'' BLAN 06-3-137237. However, the author does not support the French research policy represented by the ANR,   which promotes post-doctoral positions at the expense of permanent positions and project funding at the expense of long-term funding.
}\\
Laboratoire de math\'ematiques, UMR 8628 \\
Universit\'e Paris Sud,  Bat. 425 \\
91405 Orsay Cedex   FRANCE }

\title{Simplicity of  $\homeo(\bbD^2, \partial \bbD^2 , \mathrm{Area})$ and fragmentation of symplectic diffeomorphisms}
\date{}

\begin{document}
\maketitle

\sloppy

\begin{abstract}
In 1980, Albert Fathi asked whether the group of area-preserving homeomorphisms of the $2$-disc that are the identity near the boundary is a simple group. In this paper, we show that the simplicity of this group is equivalent to the following fragmentation property in the group of compactly supported, area preserving diffeomorphisms of the plane:
\emph{there exists a constant $m$ such that every element supported on a disc $D$ is the product of at most  $m$ elements supported on topological discs whose area are half the area of $D$.}
\end{abstract}

\selectlanguage{french}
\begin{abstract}
En 1980, Albert Fathi pose la question de la simplicité du groupe des homéomorphismes du disque qui préservent l'aire et sont l'identité près du bord. Dans cet article, nous montrons que la simplicité de ce groupe est équivalente à une propriété de fragmentation dans le groupe des difféomorphismes du plan, préservant l'aire et à support compact, à savoir :
\emph{il existe une constante $m$ telle que tout élément à support dans un disque $D$ est le produit d'au plus $m$ éléments dont les supports sont inclus dans des disques topologiques dont l'aire est la moitié de l'aire de $D$.}
\end{abstract}
\selectlanguage{english}

\textbf{AMS classification:} 37E30, 57S99, 28D15.
% 37E30 % Homeomorphisms and diffeomorphisms of planes and surfaces 
%37E45 % Rotation numbers and vectors
%37F99 % complex dyn sys.
%54H15 : transformation groups
%57S99 : Topological transformation groups (none of the above)
%28D15 General groups of measure-preserving transformations

\textbf{Keywords:} simple group; surface homeomorphism; hamiltonian dynamics; fragmentation; symplectic diffeomorphisms.

\bigskip
\bigskip
\bigskip
This paper is concerned with the algebraic study of the group
$$
G = \homeo(\bbD^2, \partial \bbD^2 , \mathrm{Area})
$$
of area-preserving homeomorphisms of the $2$-disc that are the identity near the boundary. The central open question is the following.
\begin{ques}[\cite{Fathi80}]\label{ques.fathi}
Is $G$ a simple group?
\end{ques}

The study of the simplicity of groups of homeomorphisms goes back as far as 1935. Indeed in the famous Scottish Book (\cite{ScottishBook}), S. Ulam asked if the  identity component in the group of homeomorphisms of the $n$-sphere is a simple group. This question was answered in the affirmative by Anderson and Fisher in the late fifties (\cite{Anderson58,Fisher60}).
In the seventies lots of (smooth) transformation groups were studied by D. Epstein, M. Herman, W. Thurston, J. Mather, A. Banyaga, and proved to be simple (see the books~\cite{Banyaga97} or~\cite{Bounemoura08}). 
Let us give some details on the group $G^\diff =  \diffeo(\bbD^2, \partial \bbD^2 , \mathrm{Area})$, which is the smooth analog of our group $G$. This group is not simple, since there exists a morphism from $G^\diff$ to $\bbR$, called the \emph{Calabi invariant}. But Banyaga proved that  the kernel of the Calabi invariant coincides with the subgroup $[G^\diff, G^\diff]$ generated by commutators, and is a simple group. Thus the normal subgroups of $G^\diff$ are exactly the inverse images of the subgroups of $\bbR$ under the Calabi morphism.
The analog of question~\ref{ques.fathi} is also solved in higher dimensions. Indeed A. Fathi  proved the simplicity of the group of volume preserving homeomorphisms of the $n$-ball which are the identity near the boundary, when $n \geq 3$.  However,  Question~\ref{ques.fathi} remains unsolved (see~\cite{Fathi80}).

Actually some normal subgroups of $G$ have been defined by E.~Ghys (\cite{Ghys07}, see~\cite{Bounemoura08}), and by S.~M\"uller and Y.-G.~Oh (\cite{MullerOh07}). But so far no one has been able to prove that these are proper subgroups: they might turn out to be equal to $G$.
 In this text, I propose to define still another family of normal subgroups $\{N_{\varphi}\}$ of $G$. I have not been able to prove that these subgroups are proper, but we can prove that they  are good candidates.
\begin{theo}\label{theo.good-candidates}
If some normal subgroup of the family $\{N_{\varphi}\}$ is equal to $G$, then $G$ is simple.
\end{theo}
The present work has its origin in Fathi's proof of the simplicity in higher dimensions. Fathi's argument has two steps.
The first step is a fragmentation result: any element of the group can be written as a product of two elements, each of which is supported on a topological ball whose volume is $\frac{3}{4}$ of the total volume.
The second step shows how this fragmentation property, let us call it $(P_{1})$, implies the perfectness (and simplicity) of the group.
While the second step is still valid in dimension $2$, the first one fails.
In the sequel we propose to generalise the fragmentation property $(P_{1})$  by  considering a family of fragmentation properties $(P_{\rho})$ depending on the parameter $\rho \in (0,1]$ (a precise definition is provided in section~\ref{sec.fragmentation-metric}).
A straightforward generalisation of Fathi's second step will prove that if the property $(P_{\rho})$ holds for some $\rho$, then $G$ is simple (Lemma~\ref{lem.pis} below).
On the other hand, we notice that if none of the properties $(P_{\rho})$ holds, then the subgroups $N_{\varphi}$ are proper, and thus $G$ is not simple (Lemma~\ref{lem.npins}). Thus we see firstly that Theorem~\ref{theo.good-candidates} holds, and secondly that Question~\ref{ques.fathi}  is translated into a fragmentation problem, namely the existence of some $\rho$ such that property $(P_{\rho})$ holds. 
Christian Bonatti drew my attention to the possibility of  formulating   this fragmentation problem in terms of 
a single  property $(P_{0})$. This property, which may be seen as the limit of the properties $(P_{\rho})$ as $\rho$ tends to zero, is the following:
 \emph{there exists a constant $m$ such  that any homeomorphism of the plane, supported on a disc having  area equal to one, is the composition of $m$  homeomorphisms supported on some topological discs having area equal to one half.}

This discussion is  summarised by the next theorem.

\bigskip
\bigskip

\begin{theo}\label{theo.fragmentation-homeo} The following properties are equivalent:
\begin{enumerate}
\item the group $G$ is simple,
\item there exists some $\rho \in (0,1]$ such that the property $(P_{\rho})$ holds, 
\item  property $(P_{0})$ holds.
\end{enumerate}
\end{theo}
 Furthermore we will prove that the simplicity of $G$ is also equivalent to the similar fragmentation property on the smooth subgroup  $G^\diff$ (see Lemma~\ref{lem.P0smooth} and Theorem~\ref{theo.fragmentation-diffeo} in section~\ref{sec.diffeo} below).
 We will see in section~\ref{sec.quasi-morphisms} that Entov-Polterovich quasi-morphisms, coming from Floer homology, implies that the fragmentation property $(P_{\rho})$ do not hold for $\rho \in (\frac{1}{2},1]$. 
Whether it holds or not for $\rho \in [0,\frac{1}{2}]$ remains an open question.

The definitions and precise statments are given in section~\ref{sec.fragmentation-metric}, as well as the links between properties $(P_{0})$ and $(P_{\rho})$ for $\rho>0$. The proofs of Theorem~\ref{theo.good-candidates} and~\ref{theo.fragmentation-homeo} are given in sections~\ref{sec.npins} and~\ref{sec.pis}. Sections~\ref{sec.diffeo} and~\ref{sec.profile-diffeo} provide the link with diffeomorphisms. Some more remarks, in particular the connection with  other surfaces, are mentionned in section~\ref{sec.remarks}.
Sections~\ref{sec.profile-diffeo},~\ref{sec.quasi-morphisms} and~\ref{sec.remarks} are independant.

\paragraph{Acknowledgments}
I am pleased to thank Etienne Ghys for having introduced the problem to me (in La bussière, 1997);
Albert Fathi, Yong-Geun Oh and Claude Viterbo for having organised the 2007 Snowbird conference that cast a new light on the subject;
the ``Symplexe'' team for the excellent mathematical atmosphere, and especially  Vincent Humilière, Emmanuel Opshtein and Pierre Py for the Parisian seminars and lengthy discussions around the problem; Pierre Py again for his precious commentaries on the text; Christian Bonatti for his ``je transforme ton emmental en gruyère'' trick;
and Sylvain Crovisier and Fran\c cois Béguin for the daily morning coffees, with and without normal subgroups.

%%%%%%%%%%%%%%%%%%%%%%%%
%%%%%%%%%%%%%%%%%%%%%%%%
\section{The fragmentation norms}
\label{sec.fragmentation-metric}
In the whole text, the disc $\bbD^2$ is endowed with the normalised Lebesgue measure, denoted by $\mathrm{Area}$, so that $\mathrm{Area}(\bbD^2)=1$.
 The group $G$ is endowed with the topology of uniform convergence (also called the $C^0$ topology), that turns it into a topological group. We recall that $G$ is arcwise connected: an elementary proof is provided by the famous Alexander trick (\cite{Alexander23}).
We will use the term \emph{topological disc} to denote any image of a euclidean closed disc under an element of the group $G$. As a consequence of the classical theorems by  Sch\"onflies and Oxtoby-Ulam, any Jordan curve of null area bounds a topological disc (see~\cite{OxtobyUlam41}).
 Remember that the \emph{support} of some $g \in G$ is the closure of the set of non-fixed points. For any topological disc $D$, denote by $G_{D}$ the subgroup of $G$ consisting of the elements whose support is included in the interior of $D$. Then each group $G_{D}$ is isomorphic to $G$, as shown by the following ``re-scaling'' process. Let $\Phi \in G$ be such that $D=\Phi^{-1}(D_{0})$ where $D_{0}$ is a euclidean disc. Then the map 
 $g \mapsto \Phi g \Phi^{-1}$ provides an isomorphism between the groups $G_{D_{0}}$ and $G_{D}$. We may now choose a homothecy $\Psi$ that sends  the whole disc $\bbD^2$ onto $D_{0}$, and similarly get an isomorphism $g \mapsto \Psi g \Psi ^{-1}$
  between $G$ and $G_{D_{0}}$.

%%%%%%%%%%%%%%%%%%%%%%%
\subsection*{Definition of the fragmentation metrics}
  
Let $g$ be any element of $G$. We define the \emph{size} of $g$ as follows:
$$
\mathrm{Size}(g) = \inf \{\mathrm{Area} (D), D \mbox{ is a topological disc that contains the support of } g \}.
$$
Let us emphasize the importance of the word \emph{disc}: an element $g$ which is supported on an annulus of small area surrounding a  disc of large area has a  large size. Also note that if $g$ has size less than the area of some disc $D$, then $g$ is conjugate to an element supported in $D$.

The following proposition says that the group $G$ is generated by elements of arbitrarily small size.
It is an immediate consequence of Lemma 6.5 in~\cite{Fathi80} (where the size is replaced by the diameter).
\begin{prop}[Fathi]
\label{prop.fragmentation-homeo}
Let $g \in G$, and $\rho \in (0,1]$. Then there exists some positive integer $m$,
and elements $g_{1}, \dots g_{m} \in G$ of size less than $\rho$, such that
$$
g = g_{m} \cdots g_{1}.
$$
\end{prop}
We now define the family of ``fragmentation norms''.\footnote{The definition of the fragmentation norm is not new, see example~1.24 in~\cite{BuragoIvanovPolterovich07}.}
 For any element $g \in G$ and any $\rho \in (0,1]$, we consider the least integer $m$ such that $g$ is equal to the product of $m$ elements of size less than $\rho$. This number is called the \emph{$\rho$-norm} of $g$ and is denoted by $||g||_{\rho}$.
The following properties are obvious.
\begin{prop}\label{prop.properties}
$$
||h g h^{-1}||_{\rho} = ||g||_{\rho}, \ \
||g^{-1}||_{\rho} = ||g||_{\rho}    , \ \  
||g_{1} g_{2}||_{\rho} \leq ||g_{1}||_{\rho} + ||g_{2}||_{\rho}.
$$
\end{prop}
As a consequence, the formula
$$
d_{\rho}(g_{1}, g_{2}) = ||g_{1} g_{2}^{-1}||_{\rho}
$$
defines a bi-invariant metric on $G$.

%%%%%%%%%%%%%%%%%%%%%%%
\subsection*{The normal subgroups $N_{\varphi}$}
Given some element $g \in G$, we consider the $\rho$-norm of $g$ as a function of the size~$\rho$:
$$
\rho \mapsto ||g||_{\rho},
$$
and call it the \emph{complexity profile} of $g$.  
Let $\varphi : (0,1] \to \bbR^+$ be any non-increasing function. We define the  subset  $N_{\varphi}$ containing  those elements of  $G$ whose complexity profile is essentially bounded by $\varphi$:
$$
N_{\varphi} = \{g \in G,  ||g||_{\rho} = O(\varphi(\rho))\}
$$
where the notation $\psi(\rho) = O(\varphi(\rho))$ means that there exists some $K>0$ such that $\psi(\rho) < K \varphi(\rho)$ for every small enough $\rho$.
The following is an immediate consequence of proposition~\ref{prop.properties}.
\begin{prop}
\label{prop.normal}
For any non-increasing function $\varphi : (0,1] \to \bbR^+$, the set $N_{\varphi}$ is a normal subgroup of $G$.
\end{prop}

The reader who wants some examples where we can estimate the complexity profile may jump to section~\ref{sec.profile-diffeo}, where we will see that commutators of diffeomorphisms have a profile equivalent to the function $\varphi_{0} : \rho \mapsto \rho^{-1}$. This will imply that $N_{\varphi_{0}}$ is the smallest non-trivial subgroup of our family $\{N_{\varphi}\}$.

%%%%%%%%%%%%%%%%%%%%%%%
\subsection*{The fragmentation properties $(P_{\rho})$}

Let $\rho \in (0,1]$. We now define our fragmentation property $(P_{\rho})$  by asking for a uniform bound in the fragmentation of elements of size less than $\rho$ into elements of a smaller given size.
\begin{itemize}
\item[$\mathbf{(P_{\rho})}$]  There exists some number $s \in (0,\rho)$, and some positive integer $m$, such that any $g\in G$ of size less than $\rho$ satisfies $||g||_{s} \leq m$.
\end{itemize}

Here are some easy remarks. Let us denote by $P(\rho,s)$ the property that there exists a bound $m$ with  $||g||_{s} \leq m$ for every element $g$ of size less than $\rho$. Fix some $\rho \in (0,1]$ and some ratio $k \in (0,1)$.
Assume that property $P(\rho,k\rho)$ holds. Then  by re-scaling we get  that  property $P(\rho',k\rho')$ also holds for any $\rho' < \rho$ (with the same bound $m$).
In particular we can iterate the fragmentation to get, for every positive $n$,
property $P(\rho, k^n\rho)$ (with the bound  $m^n$).
This shows that  property $P(\rho,s)$ implies property 
 $P(\rho,s')$ for every $s' < s$. The converse is clearly true, so that property $P(\rho,s)$ depends only on $\rho$ and not on $s$.
In particular we see that  property $(P_{\rho})$ is equivalent to the existence of a number $m$ such that every $g$ of size less than $\rho$ satisfies
$$
||g||_{\frac{\rho}{2}} \leq m.
$$
Also note that property $P(\rho_{0})$ implies property $P(\rho_{1})$ if $\rho_{1} < \rho_{0}$ (again by re-scaling). Thus property $P_{\rho}$ is more and more likely to hold as $\rho$ decreases from $1$ to $0$.

%%%%%%%%%%%%%%%%%%%%%%%
\subsection*{The fragmentation property $(P_{0})$}

In this paragraph we introduce the fragmentation property $(P_{0})$, and prove that it is equivalent to the existence of some $\rho>0$ such that property $(P_{\rho})$ holds.
Consider, just for the duration of this section, the group
$$
\homeo_{c}(\bbR^2, \mathrm{Area})
$$
of compactly supported, area preserving homeomorphisms of the plane.
Any image of a euclidean closed disc under some  element of this bigger group will again be called a topological disc.
We also define the size of an element of the group  as in $G$. Property $(P_{0})$ is as follows.
\begin{itemize}
\item[$\mathbf{(P_{0})}$]  There exists some positive integer $m$ such that any $g\in \homeo_{c}(\bbR^2, \mathrm{Area})$ of size less than $1$ 
is the composition of at most $m$ elements of $\homeo_{c}(\bbR^2, \mathrm{Area})$ of size less than $\frac{1}{2}$.
\end{itemize}
%Contrarily to properties $(P_{\rho})$, here we assume no bound on the area spanned by the supports of the pieces of the fragmentation.
Since each piece of the fragmentation provided by property $(P_{0})$ is supported on a disc with area $\frac{1}{2}$, the union of the supports has area at most $\frac{m}{2}$, but this gives no bound on the area of a topological disc containing this union. However, if the union of the supports surround a region with big area, then we may find a new fragmentation 
by ``bursting the bubble'', \emph{i.e.} conjugating the situation by a map that contracts the areas of the surrounded regions and preserves the area everywhere else. This is the key observation, due to Christian Bonatti,  to the following lemma.
\begin{lemm}\label{lem.P0}
Property $(P_{0})$ holds if and only if there exists some $\rho \in (0,1]$ such that property $(P_{\rho})$ holds.
\end{lemm}

\begin{proof}
Let us prove the easy part. Suppose $(P_{\rho})$ holds for some $\rho>0$, let $m$ be a bound for $||g||_{\frac{\rho}{2}}$ for those $g \in G$ of size less than $\rho$. Let $g \in \homeo_{c}(\bbR^2, \mathrm{Area})$ of size less than $1$.
Choose some element of the group that sends the support of $g$ into the euclidean unit disc $\bbD^2$, and compose it with the homothecy that sends $\bbD^2$ onto a disc of area $\rho$ included in $\bbD^2$; we denote by $\Psi$ the resulting map.
 Then $\Psi g \Psi^{-1}$ is an element of $G$ of size less than $\rho$. According to hypothesis $(P_\rho)$,
 we may write this element as a composition of  $m$ elements of $G$ of size less than ${\frac{\rho}{2}}$.
We may conjugate these elements by $\Psi^{-1}$ and take the composition to get a fragmentation of $g$ into $m$ elements of size $\frac{1}{2}$. Thus $(P_{0})$ holds.

Now assume that $(P_{0})$ holds, and let $m \geq 2$ be  given by this property. We will prove that property $(P_{\rho})$ holds for $\rho = \frac{2}{m}$.
Consider, in the plane, a euclidean disc $D$ of area $\frac{1}{\rho}$.
By the same re-scaling trick as before, it suffices to prove that any $g \in \homeo_{c}(\bbR^2, \mathrm{Area})$ with size less than $1$ and supported in the interior of $D$ may be fragmented as a product of $m$ elements of $\homeo_{c}(\bbR^2, \mathrm{Area})$ with size less than one half and supported in the interior of $D$.
Property $(P_{0})$ provides us with a fragmentation $g = g'_{m} \circ \cdots \circ g'_{1}$
by elements of size less than one half, but maybe not supported in $D$. Now comes the ``bursting the bubbles'' trick.
Let $D'$ be a topological disc whose interior contains all the supports of the $g'_{i}$'s.
The union of the supports has area less than $\frac{1}{\rho}$. Thus we may find some topological discs 
$K'_{1}, \dots K'_{\ell}$, included in the interior of $D'$, that are pairwise disjoint and disjoint from the supports of the $g'_{i}$'s, such that
$$
\mathrm{Area}\left (D' \setminus \bigcup_{j=1}^\ell K'_{j} \right) < \frac{1}{\rho}.
$$
Denote by $D_{0}$ the support of our original map $g$. Note that $D_{0}$ is included in the union of the supports of the $g'_{i}$'s, thus it is disjoint from the $K'_{j}$'s.
Since $D$ has area $\frac{1}{\rho}$, the previous inequality ensures the existence of some pairwise disjoint discs $K_{1}, \dots K_{\ell}$ in the interior of $D$, disjoint from $D_{0}$,  such that 
$$
\mathrm{Area}\left (D \setminus \bigcup_{j=1}^\ell K_{j} \right)  = \mathrm{Area}\left (D' \setminus \bigcup_{j=1}^\ell K'_{j} \right).
$$
Using Sch\"onflies and Oxtoby-Ulam theorems, we can construct a homeomorphism $\Psi$ of the plane satisfying the following properties:

\begin{enumerate}
\item $\Psi$ is the identity on $D_{0}$,
\item $\Psi(D') = D$, $\Psi(K'_{j})= K_{j}$ for each $j$,
\item the restriction of $\Psi$ to the set $D' \setminus \cup_{j=1}^\ell K'_{j}$ preserves the area.
\footnote{\label{foot.smooth} Having in mind the smooth case (Lemma~\ref{lem.P0smooth} below), we notice that we  may further demand that the map $\Psi$ is a $C^\infty$-diffeomorphism on the interior of $D' \setminus \cup_{j=1}^\ell K'_{j}$. Actually, we may even choose the sets $D'$, $K_{j}$ and $K'_{j}$ to be smooth discs, and then the map $\Psi$ may be chosen to be a  $C^\infty$-diffeomorphism of the plane.}
\end{enumerate}

The first item shows that $\Psi g \Psi^{-1} = g$. Now for each $i$ we define $g_{i} = \Psi g'_{i} \Psi^{-1}$. Then the second item guarantees that the $g_{i}$'s are supported in $D$, and the third item entails that they preserve area and have size less than one half. The product of the $g_{i}$'s is equal to $g$, which provides the desired fragmentation.
\end{proof}

%%%%%%%%%%%%%%%%%%%%%%%%
%%%%%%%%%%%%%%%%%%%%%%%%
\section{Simplicity implies fragmentation}
\label{sec.npins}
\begin{lemm}\label{lem.npins} 
Assume that none of the properties $(P_{\rho}),\rho \in (0,1]$ holds. 
Let $\varphi : (0,1] \to \bbR^+$ be any function. Then the normal subgroup $N_{\varphi}$ is proper, \emph{i. e.} it is not equal to $G$. In this case
the group $G$ is not simple.
\end{lemm}
If we consider any element $f \neq \mathrm{Id}$ in $G$, and  the function $\varphi_{f} : \rho \mapsto  ||f||_{\rho}$, then   the normal subgroup $N_{\varphi_{f}}$ contains $f$ and thus is not equal to $\{\mathrm{Id}\}$. Hence the non-simplicity of $G$ will be a consequence of the non-triviality of the subgroups $N_{\varphi}$.

\begin{proof}
According to the easy remarks following the definition of property $(P_{\rho})$, the hypothesis of the lemma reads the following way:
\begin{itemize}
\item[$(\star)$] for every $\rho \in (0,1]$ and every positive integer $m$ there exists some element $g$ of size less than $\rho$ such that $ ||g||_{\frac{\rho}{2}} > m$.
\end{itemize}

 We fix any function  $\varphi : (0,1] \to \bbR^+$, and we will construct some element $g$ in $G$ that does not belong to $N_{\varphi}$. 
Let us define $D_{0}= \bbD^2$. We  pick two sequences of discs $(C_{i})_{i \geq 1}$ and $(D_{i})_{i \geq 1}$ converging to a point, such that for every~$i$ (see figure~\ref{fig.sequence1}),
\begin{itemize}
\item[--] $C_{i}$ and  $D_{i}$ are disjoint and included in $D_{i-1}$, 
\item[--] the area of $D_{i}$ is less than half the area of $C_{i}$.
\end{itemize}
\begin{figure}[htbp]
\begin{center}
\def\JPicScale{1}
\includegraphics{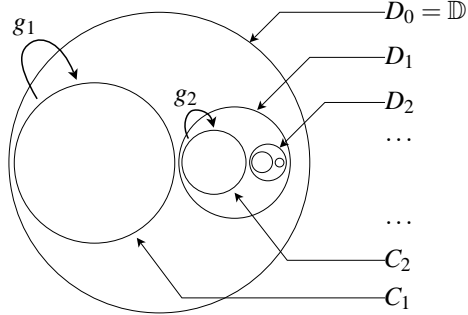}

\caption{Construction of $g$}
\label{fig.sequence1}
\end{center}
\end{figure}

We denote the area of $C_{i}$  by $\rho_{i}$.
We will construct a sequence $(g_{i})_{i \geq 1}$, with each $g_{i}$ supported in the interior of $C_{i}$, and then  $g$ will be defined as 
the (infinite)  product of the $g_{i}$'s. Note that since the discs $C_{i}$'s are pairwise disjoint this product has a meaning, and since the sequence $(C_{i})$ converges to a point it actually defines an element of $G$.
Since all the $g_{j}$'s with $j>i$ will be supported in the interior of the disc $D_{i}$ whose area is less than $\rho_{i}/2$ we will get
\begin{eqnarray}
||g||_{\frac{\rho_{i}}{2}} & \geq & ||g_{i} \dots g_{1}||_{\frac{\rho_{i}}{2}}-1.
\end{eqnarray}
The sequence $(g_{i})$ is constructed by induction. Assume $g_{1}, \dots , g_{i-1}$ have been constructed. Using hypothesis $(\star)$, we may choose $g_{i}$ supported on $C_{i}$ such that $
||g_{i}||_{\frac{\rho_{i}}{2}}$ is arbitrarily high, more precisely we demand the following inequality:
\begin{eqnarray}
||g_{i}||_{\frac{\rho_{i}}{2}} & \geq &  \frac{1}{\rho_{i}} \varphi\left(\frac{\rho_{i}}{2}\right) \quad + \quad  
||g_{i-1} \dots g_{1}||_{\frac{\rho_{i}}{2}}
\quad + \quad  1.
\end{eqnarray}
Using inequality (1), the triangular inequality and inequality (2) we get
\begin{eqnarray*}
||g||_{\frac{\rho_{i}}{2}} &\geq &   ||g_{i} \dots g_{1}||_{\frac{\rho_{i}}{2}}-1 \\
 & \geq & ||g_{i}||_{\frac{\rho_{i}}{2}}  
 \quad - \quad  
  ||g_{i-1} \dots g_{1}||_{\frac{\rho_{i}}{2}} -1 \\
 &   \geq &  \frac{1}{\rho_{i}} \varphi\left(\frac{\rho_{i}}{2}\right).
\end{eqnarray*}
This proves that the complexity profile of $g$ is not equal to $O(\varphi)$.
In other words $g$ does not belong to $N_{\varphi}$.
\end{proof}

%%%%%%%%%%%%%%%%%%%%%%%%
%%%%%%%%%%%%%%%%%%%%%%%%
\section{Fragmentation implies simplicity}
\label{sec.pis}

\begin{lemm}\label{lem.pis}
Assume that property $(P_{\rho})$ holds for some $\rho \in (0,1]$.
Then $G$ is simple.
\end{lemm}
This lemma is just a slight generalisation of Fathi's argument showing that, under property $(P_{1})$, the group  $G$ is perfect: any element decomposes as a product of commutators.  Then  perfectness implies simplicity: this is due to ``Thurston's trick'',  for completeness the argument is included in the proof below.

\begin{proof}
We assume that there exists a number $\rho \in (0,1]$ and a positive integer $m$
such that any element of size less than $\rho$ may be written as the product of $m$ elements of size less than $\frac{\rho}{2}$.

 Let $C_{1}$ be a small disc. By usual fragmentation (proposition~\ref{prop.fragmentation-homeo}), any element of $G$ is a product of elements supported in a disc of area less than that of $C_{1}$, and any such element is conjugate to an element supported in the interior of $C_{1}$. Thus to prove perfectness it is enough to consider some element $g$ supported in the interior of $C_{1}$ and to prove that $g$ is a product of commutators.

Let us first prove that such a $g$ is a product of two commutators when considered in the group $ \homeo(\bbD^2, \partial \bbD^2 )$, that is, let us forget for a while about the area (this is  a ``pedagogical'' step).
Choose two sequences of   discs $(C_{i})_{i \geq 1}$ and $(D_{i})_{i \geq 1}$ converging to a point, such that (see figure~\ref{fig.sequence2})
\begin{itemize}
\item[--]  the interior of $D_{i}$ contains both $C_{i}$ and $C_{i+1}$, 
\item[--] the $C_{i}$'s  are pairwise disjoint,
\item[--] the $D_{2i}$'s (resp. the $D_{2i+1}$'s) are pairwise disjoint.
\end{itemize}
\begin{figure}[htbp]
\begin{center}
\def\JPicScale{1}
\includegraphics{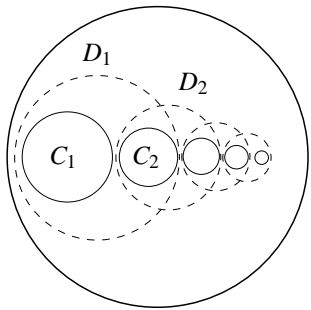}
\caption{The sequences $(C_{i})_{i \geq 1}$ and $(D_{i})_{i \geq 1}$}
\label{fig.sequence2}
\end{center}
\end{figure}
 For any $i\geq 1 $ choose some $h_{i} \in\homeo(\bbD^2, \partial \bbD^2 )$, supported on $D_{i}$, that sends $C_{i}$ onto $C_{i+1}$. We let $g_{1}:= g$, thus $g_{1}$ is supported on $C_{1}$, and define inductively $g_{i+1} := h_{i} g_{i} h_{i }^{-1}$~; thus $g_{i}$ is a ``copy'' of $g$, supported on $C_{i}$, and the $g_{i}$'s are pairwise commuting.
Let 
$$
K := g_{2} g_{3}^{-1} g_{4} g_{5}^{-1} \cdots , \ \  K' := g_{1} g_{2}^{-1} g_{3} g_{4} ^{-1} \cdots
$$
so that $KK'=K'K=g$. The map  $K = [g_{2},h_{2} ][ g_{4}, h_{4}] \cdots $ may be seen as an infinite product of commutators, but we need a finite product.
Now define 
$$
G :=  g_{2} g_{4} \dots , \ \ H:= h_{2} h_{4}  \cdots , \ \   G' := g_{1} g_{3} \cdots , \ \  H' := h_{1} h_{3} \cdots 
$$
and observe that $K = [G,H]$ and $K'=[G',H']$: indeed these equalities may be checked independently on each disc $D_{i}$.
Thus $g = [G,H][G',H']$ is a product of two commutators in $\homeo(\bbD^2, \partial \bbD^2 )$.

Now let us take care of the area. We will use sequences $(C_{i})$ and $(D_{i})$ as before, and we will get around the impossibility of shrinking $C_{i}$ onto $C_{i+1}$ inside the group $G$ by using the fragmentation hypothesis.
We may assume, for every $i$, the equality
$$
\area(C_{i+1}) = \frac{1}{2} \area(C_{i}).
$$
Moreover, by fragmentation, we may assume this time that $g$ is supported in the interior of  a disc $C'_{1} \subset C_{1}$ of area  $\rho\area (C_{1})$.
We use the fragmentation hypothesis re-scaled on $C_{1}$ to write 
$$
g = f_{1,1} \dots f_{1,m}
$$
(see figure~\ref{fig.perfect}) with each $f_{1,j}$ supported in the interior of a topological disc included in $C_{1}$ and whose area is 
$$
\frac{\rho}{2}\area (C_{1}) = \rho \area (C_{2}).
$$
We choose a disc $C'_{2} \subset C_{2}$ whose area also equals $\rho \area (C_{2})$ and, for each $j=1, \dots , m$, some $h_{1,j} \in G$ supported on $D_{1}$ and sending the support of $f_{1,j}$ inside $C'_{2}$. 
We define 
$$
g_{2} := \prod_{j=1, \dots , m} h_{1,j} f_{1,j} h_{1,j}^{-1}
$$
which is supported on $C'_{2}$ (see figure~\ref{fig.perfect}).
\begin{figure}[htbp]
\begin{center}
\def\JPicScale{1.4}
\includegraphics{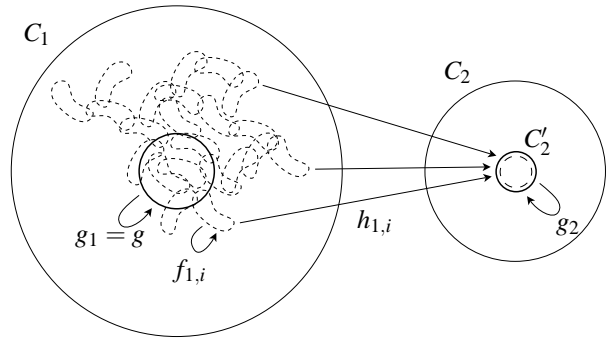}
\caption{Fragment and push every piece inside the small disc...}
\label{fig.perfect}
\end{center}
\end{figure}
 We apply recursively the (re-scaled) fragmentation hypothesis to get a sequence $(g_{i})_{i \geq 1}$ with each $g_{i}$ supported in the interior of a  disc $C'_{i} \subset C_{i}$ having area $\rho \area(C_{i})$ and sequences $(f_{i,j})_{i \geq 1,j=1, \dots , m}$ and $(h_{i,j})_{i \geq 1,j=1, \dots , m}$ with $f_{i,j}$ supported on $C_{i}$ and $h_{i,j}$   supported on $D_{i}$, such that
$$
g_{i} = \prod_{j=1, \dots , m}  f_{i,j}  \mbox{ and } 
g_{i+1} = \prod_{j=1, \dots , m} h_{i,j} f_{i,j} h_{i,j}^{-1}.
$$
Obviously $g_{i}$ and $g_{i+1}$ are equal up to a product of commutators \emph{whose number of terms  depends only on $m$}. More precisely, we have 
\begin{eqnarray*}
g_{i} g_{i+1} ^{-1} & =  & \prod_{j=1, \dots , m}  f_{i,j} \prod_{j=m, \dots , 1} h_{i,j} f_{i,j}^{-1} h_{i,j}^{-1} \\
  & = & \left[ f_{i,1} , P \right]  %h_{i,1} f_{i,1}^{-1} h_{i,1}^{-1} f_{i,1} 
\left( \prod_{j=2, \dots , m}  f_{i,j}  \prod_{j=m, \dots , 2} h_{i,j} f_{i,j} ^{-1} h_{i,j}^{-1} \right)
 \left[f_{i,1}, h_{i,1} \right]
\end{eqnarray*}
where $P$ is equal to the term between parentheses,
 and we see recursively that $g_{i} g_{i+1} ^{-1}$ is a product of $2m$ commutators
of elements  supported in $D_{i}$;  we write
 $$
 g_{i} g_{i+1} ^{-1} = \prod_{j=1, \dots , 2m} [s_{i,j}, t_{i, j}].
 $$
It remains to define the infinite commutative products 
$$
K := g_{2} g_{3}^{-1} g_{4} g_{5}^{-1} \cdots , \ \  K' := g_{1} g_{2}^{-1} g_{3} g_{4} ^{-1} \cdots
$$
$$
S_{j} :=  s_{2,j} s_{4,j} \dots , \ \ T_{j}:= t_{2,j} t_{4,j}  \cdots , \ \   S'_{j} := s_{1,j} s_{3,j} \cdots , \ \  T'_{j} := t_{1,j} t_{3,j} \cdots 
$$
and to check that 
$$
K= \prod_{j = 1, \dots 2m} [S_{j}, T_{j}], \quad
K' = \prod_{j = 1, \dots 2m} [S'_{j}, T'_{j}], \quad
\mbox{ and } g = KK' 
$$
is a product of $4m$ commutators.
This proves that  $G$ is perfect.

Let us recall briefly, according to Thurston,  how perfectness implies simplicity.
Let $D$ be a disc and $g,h\in G$ be such that the discs $D, g(D), h(D)$ are pairwise disjoint.
Let $u,v\in G$ be supported in $D$. In this situation the identity
$$
[u,v] = [[u,g],[v,h]]
$$
may easily be checked, and shows that $[u,v]$ belongs to the normal subgroup generated by $g$.
Now given any $g\neq \mathrm{Id}$ in $G$, one can find an $h \in G$ and a disc $D$ such that the above situation takes place.
If $G$ is perfect then so is the isomorphic group $G_{D}$, hence every $f$ supported in $D$ is a product of commutators supported in $D$, and by the above equality such an $f$ belongs to the normal subgroup generated by $g$. By fragmentation this subgroup is thus equal to $G$. This proves that $G$ is simple, and completes the proof of the lemma.
\end{proof}

%%%%%%%%%%%%%%%%%%%%%%%%
%%%%%%%%%%%%%%%%%%%%%%%%
\section{Fragmentation of diffeomorphisms}
\label{sec.diffeo}
Here we further translate Question~\ref{ques.fathi}  into the  diffeomorphisms subgroup $G^\diff$. 

Let $G^\diff = \diffeo(\bbD^2, \partial \bbD^2 , \mathrm{Area})$ be the group of elements of $G$ that are $C^\infty$-diffeomorphisms. 
Note that for every topological disc $D$  the group of elements supported in the interior of $D$,
$$G_{D}^\diff := G^\diff \cap G_{D}, $$
 is isomorphic to $G^\diff$, even when $D$ is not smooth: indeed we may use Riemann conformal mapping theorem and Moser's lemma to find a smooth diffeomorphism $\Phi$ between the interiors of $\bbD^2$ and $D$ with constant Jacobian, and the conjugacy by $\Phi$ provides an isomorphism (see~\cite{GreeneShiohama79} for the non-compact version of Moser's lemma).
As in the continuous case, we define the $\rho$-norm $||g||^\diff_{\rho}$ of any element $g \in G^\diff$
as the minimum number $m$ of elements $g_{1}, \dots , g_{m}$ of $G^\diff$, having size less than $\rho$, whose composition is equal to $g$.
The fragmentation properties $(P^\diff_{\rho})$ are defined as in the continuous setting.
\begin{itemize}
\item[$\mathbf{(P^\diff_{\rho})}$] (for $\rho \in (0,1]$)  There exists some number $s \in (0,\rho)$, and some positive integer $m$, such that any $g\in G^\diff$ of size less than $\rho$ satisfies $||g||^\diff_{s} \leq m$.
\item[$\mathbf{(P^\diff_{0})}$]  There exists some positive integer $m$ such that any $g\in \diffeo_{c}(\bbR^2, \mathrm{Area})$ of size less than $1$  is the composition of at most $m$ elements of $\diffeo_{c}(\bbR^2, \mathrm{Area})$ of size less than $\frac{1}{2}$.
\end{itemize}
Here $\diffeo_{c}(\bbR^2, \mathrm{Area})$ is the group  of  $C^\infty$-diffeomorphisms of the plane that are compactly supported and preserve the area. The smooth version of Lemma~\ref{lem.P0} holds, with the same proof (using footnote~\ref{foot.smooth}).
\begin{lemm}\label{lem.P0smooth}
Property $(P^\diff_{0})$ holds if and only if there exists some $\rho \in (0,1]$ such that property $(P^\diff_{\rho})$ holds.
\end{lemm}

We now turn to the equivalence between fragmentation properties for homeomorphisms and diffeomorphisms.
\begin{theo}\label{theo.fragmentation-diffeo}
For any $\rho \in (0,1]$, the properties $(P_{\rho})$ and $(P^\diff_{\rho})$ are equivalent.
\end{theo}
The proof of this equivalence requires two ingredients. The first one is the  density of $G^\diff$ in $G$; this is a classical result, see for example~\cite{Sikorav07}. The second one is the uniformity of the fragmentation in $G$, and in $G^\diff$, inside some $C^0$-neighbourhood of the identity. This is provided by the following proposition, which is proved below.
\begin{prop}
\label{prop.fragmentation-locale}
For any $\rho \in (0,  1)$, there exists a neighbourhood $\cV_{\rho}$ of the identity in $G$ with the following properties.

\begin{itemize}
\item[--] Any $g$ in $\cV_{\rho}$  satisfies $||g||_{\rho} \leq \frac{2}{\rho}$;
\item[--] any $g$ in $\cV_{\rho} \cap G^\diff$  satisfies $||g||^\diff_{\rho} \leq \frac{2}{\rho}$.
\end{itemize}
\end{prop}

As a consequence of this proposition we get the following comparison between the norms $||.||^\diff_{\rho}$ and $||.||_{\rho}$ on $G^\diff$.
\begin{coro}\label{coro.comparison-metrics}
Any $g \in G^\diff$ satisfies
$$
||g||_{\rho} \leq ||g||^\diff_{\rho} \leq ||g||_{\rho} + \frac{2}{\rho}.
$$
\end{coro}

\begin{proof}[Proof of the corollary]
 The first inequality is clear. To prove the second one consider some $g \in G^\diff$, and let $m = ||g||_{\rho}$.
 By definition there exists some elements $g_{1}, \dots , g_{m}$ in $G$ of size less than $\rho$ such that 
 $g = g_{m} \cdots  g_{1}$. Since  $G^\diff$ is $C^0$-dense in $G$, for any topological disc $D$ the subgroup
 $G^\diff_{D}$ is also $C^0$-dense in $G_{D}$.
 Thus we can find elements $g'_{1}, \dots , g'_{m}$ in $G^\diff$, still having  size less than $\rho$, whose product
 $g' = g'_{m} \cdots  g'_{1}$ is  a diffeomorphism arbitrarily $C^0$-close to $g$. According to the second item in proposition~\ref{prop.fragmentation-locale}
we get  
$$
||g'g^{-1}||^\diff_{\rho} \leq \frac{2}{\rho}.
$$
Since $||g'||^\diff_{\rho} \leq m$ the triangular inequality gives $||g||^\diff_{\rho} \leq m + \frac{2}{\rho}$, as wanted.
\end{proof}

\begin{proof}[Proof of Theorem~\ref{theo.fragmentation-diffeo}]
Let us fix $\rho \in (0,1]$.
The fact that $(P_{\rho})$ implies $(P_{\rho}^\diff)$ immediatly follows from the corollary.
We prove the converse implication.  Assume that $(P_{\rho}^\diff)$ holds: there exists some $s \in (0,\rho)$ and $m>0$ such that
the quantity $||g'||^\diff_{s}$ is bounded by $m$ on the elements $g'$ of $G^\diff$ of size less than $\rho$.
Choose any $g \in G$ of size less than $\rho$.
According to the first item in proposition~\ref{prop.fragmentation-locale}, and using the density of $G^\diff$ in $G$, 
we get some $g' \in G^\diff$ having size less than $\rho$, and sufficiently close to $g$ so that 
$||g'g^{-1}||_{s} \leq \frac{2}{s}$. Since $(P_{\rho}^\diff)$ holds we also have $||g'||^\diff_{s} \leq m$, and thus
$||g'||_{s} \leq m$, from which we get $||g||_{s} \leq m + \frac{2}{s}$. Thus $(P_{\rho})$ also holds.
\end{proof}

We now turn to the proof of proposition~\ref{prop.fragmentation-locale}.
The classical proof of fragmentation for diffeomorphisms relies on the inverse mapping theorem and would only gives uniformity in a $C^1$-neighbourhood of the identity (see~\cite{Banyaga97,Bounemoura08}). Thus we will rather try to mimic the proof of the fragmentation for homeomorphisms.

\begin{proof}[Proof of proposition~\ref{prop.fragmentation-locale}]
The proof of both items (uniformity of local fragmentation for homeomorphisms and diffeomorphisms) are very similar; we will only provide details for the diffeomorphisms case. 

We choose an integer $m$ bigger than $\frac{2}{\rho}$, and we cut the disc into $m$ strips of area less than $\frac{\rho}{2} $: more precisely, we choose $m$ topological discs $D_{1}, \cdots, D_{m}$ such that (see figure~\ref{fig.Vrho})
 \begin{enumerate}
\item $\mathrm{Area}(D_{i}) \leq \frac{\rho}{2}$,
\item $\bbD^2 = D_{1} \cup \cdots \cup D_{m}$,
\item $D_{i} \cap D_{j} = \emptyset $ if $| j-i| > 1$,
\item $D_{i} \cup D_{i+1} \cup \cdots \cup D_{j}$  is a topological disc for every $i \leq j$, and the intersection  of $D_{1} \cup \cdots \cup D_{i}$  and $D_{i} \cup \cdots \cup D_{m}$ with the boundary of $\bbD^2$ is non-empty and connected for every $i$.
\end{enumerate}
Now we define the following set $\cV_{\rho}$:
$$
\cV_{\rho} = \left\{g \in G   \mbox{ such that } g(D_{i}) \cap D_{j} = \emptyset   \mbox{ for every } i,j \mbox{ with } | j-i| > 1\right\}.
$$
\begin{figure}[htbp]
\begin{center}
\def\JPicScale{1}
\includegraphics{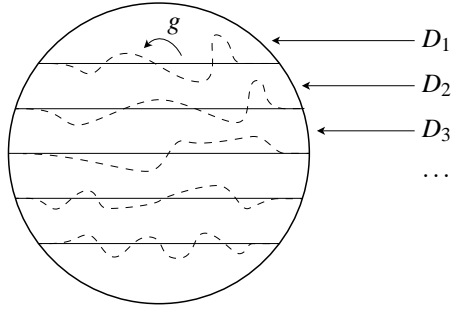}

\caption{The discs $D_{i}$ and the action of some $g$ in $\cV_{\rho}$}
\label{fig.Vrho}
\end{center}
\end{figure}

 Note that, due to item~3,  $\cV_{\rho}$ is a $C^0$-neighbourhood of the identity in $G$. We will prove that each element of $\cV_{\rho} \cap G^\diff$ can be written as a product of $m-1$ elements of $G_{\rho}^\diff$.

Let $g \in \cV_{\rho} \cap G^\diff$. By hypothesis  $D_{1}$ and  $g(D_{1})$  are both  disjoint from the topological disc $D_{3} \cup \cdots \cup D_{m}$. By the classical Lemma~\ref{lem.classical} below, we can find $\Psi_{1} \in G^\diff$ such that 
\begin{itemize}
\item[--]$\Psi_{1} = g$ on some neighbourhood of $D_{1}$,
\item[--]$\Psi_{1}$ is the identity on some neighbourhood of  $D_{3} \cup \cdots \cup D_{m}$.
\end{itemize}
The diffeomorphism $\Psi_{1}$ is supported in the interior of the topological disc $D_{1} \cup D_{2}$ whose area is less than or equal to $\rho$.
Let $g_{1} := \Psi_{1}^{-1} g$, thus $g_{1}$ is supported in the interior of $D_{2} \cup \cdots \cup D_{m}$, and we easily check that this diffeomorphism is still in $\cV_{\rho}$. In particular $D_{1} \cup D_{2}$ and its image under $g_{1}$ are both disjoint from $D_{4} \cup \cdots \cup D_{m}$. We apply again the lemma to get some $\Psi_{2}\in G$ such that 
\begin{itemize}
\item[--]$\Psi_{2} = g_{1}$ on some neighbourhood of $D_{1}\cup D_{2}$,
\item[--]$\Psi_{2}$ is the identity on some neighbourhood of  $D_{4} \cup \cdots \cup D_{m}$.
\end{itemize}
Thus $\Psi_{2}$ is supported in the interior of $D_{2} \cup D_{3}$. Let 
$g_{2} := \Psi_{2}^{-1} g_{1}$; this diffeomorphism is in $\cV_{\rho}$ and is supported in the interior of $D_{3} \cup \cdots \cup D_{m}$.
In the same way we construct diffeomorphisms $\Psi_{1}, \dots ,\Psi_{m-1}$, 
such that each $\Psi_{i}$ is supported in the interior of $D_{i} \cup D_{i+1}$, and such that 
$g = \Psi_{1} \circ \cdots \circ \Psi_{m-1}$. This completes the proof for the diffeomorphisms case.
\end{proof}

%%%%%%%%%%
\begin{lemm}\label{lem.classical}
Let $D'_{1}, D'_{2}$ be two disjoint  topological discs in $\bbD^2$, and assume that the intersection of $D'_{1}$ (resp. $D'_{2}$)  with the boundary of $\partial \bbD^2$ is non-empty and connected (and thus $\bbD^2 \setminus \inte(D'_{1} \cup D'_{2})$ is again a topological disc).
Let $\Phi \in G^\diff$, and suppose that  $\Phi(D'_{1})$ is disjoint from $D'_{2}$.
Then there exists $\Psi \in G^\diff$  such that $\Psi = \Phi$ on some neighbourhood of  $D'_{1}$ and $\Psi = \mathrm{Id}$ on some neighbourhood of $D'_{2}$.
\end{lemm}
\begin{proof}[Proof of the lemma]
By Smale's theorem (\cite{Smale59}) and Moser's lemma (see for example~\cite{Banyaga97} or~\cite{Bounemoura08}), the group $G$ is arcwise connected. Let $(\Phi_t)_{t \in [0,1]}$ be a smooth isotopy from the identity to $\Phi$ in $G$.
It is easy to find another smooth isotopy  $(g_t)_{t \in [0,1]}$ (that does not preserve area), supported in the interior of $\bbD^2$,
such that    
\begin{itemize}
\item[--]for every $t$,  $g_t (\Phi_t (D'_{1}))$ is disjoint from $D'_{2}$,
\item[--]$g_0 = g_1= \mathrm{Id}$.
\end{itemize}

The isotopy $(\Phi'_{t}) = (g_t \Phi_t)$ still goes from  the identity to $\Phi$.
Consider the vector field tangent to this isotopy, and multiply it by some smooth function 
that is equal to $1$ on some neighbourhood of $\cup_{t} \Phi'_{t}(D'_{1})$ and vanishes on $D'_{2}$.
By integrating this truncated vector field we get  another isotopy 
 $(\Psi_t)$ such that 
\begin{itemize}
\item[--] on some neighbourhood of $D'_{1}$ we have $\Psi_t = g_t \Phi_t$ for every $t$, and in particular $\Psi_1 = \Phi$,
\item[--]the support of $\Psi_{1}$ is disjoint from $D'_{2}$.
\end{itemize}
Thus $\Psi_{1}$ satisfies the conclusion of the lemma, except that it does not preserve the area. Let $\omega_{0}$ be the Area form on $\bbD^2$, and $\omega_{1}$ be the pre-image of  $\omega_{0}$ under $\Psi_{1}$.
Then $\omega_{1} = \omega_{0}$ on some neighbourhood of $\partial \bbD^2 \cup D'_{1} \cup D'_{2}$. By Moser's lemma we may find some $\Psi_{2}\in G$, whose support is disjoint from $D'_{1}$ and $D'_{2}$, and that sends $\omega_{1}$ to $\omega_{0}$. The diffeomorphism 
$\Psi = \Psi_{1} \Psi_{2}^{-1}$ suits our needs.
\end{proof}

Note that this lemma has a $C^0$-version, which is proved by replacing Smale's theorem by Alexander's trick, the truncation of vector fields by Sch\"onflies's theorem, and Moser's lemma by Oxtoby-Ulam's theorem.

%%%%%%%%%%%%%%%%%%%%%%%%
%%%%%%%%%%%%%%%%%%%%%%%%
\section{Profiles of diffeomorphisms}\label{sec.profile-diffeo}
One can wonder what the complexity profile looks like for a diffeomorphism, both inside the group $G^\diff$ and inside the group $G$. 
The following proposition only partially solves this problem.

\begin{prop}~
\label{prop.profile-diffeo}
\begin{itemize}
\item[--] For any $g \in G^\diff$ we have 
$||g||^\diff_{\rho} =O(\frac{1}{\rho^2})$.
\item[--] For any $g$ in the commutator subgroup $[G^\diff, G^\diff]$ we have 
$||g||^\diff_{\rho} =O(\frac{1}{\rho})$.
\end{itemize}
\end{prop}
If the support of $g \in G$ has area $A$, then clearly for  any $\rho$ we need at least $\frac{A}{\rho}$ elements of $G_{\rho}$ or $G^\diff_{\rho}$ to get $g$.
Thus, according to the second point, the profile of any $g$ in $[G^\diff, G^\diff]$ is 
bounded from above and below by multiples of the function $\varphi_{0} : \rho \mapsto \frac{1}{\rho}$; and this holds both in $G$ and $G^\diff$. In particular we see that $N_{\varphi_{0}}$ is the smallest non-trivial subgroup of our family $\{N_{\varphi}\}$.
  I have not been able to decide whether  the first point is optimal, nor whether  $G^\diff \subset N_{\varphi_{0}}$ or not (if not, of course, then $G$ is not simple).
  It might also happens that $G^\diff$ is included in  $N_{\varphi_{0}}$ but not in the analog smooth group $N^\diff_{\varphi_{0}}$.

We also notice that \emph{every non trivial normal subgroup of $G$ contains the commutator subgroup $[G^\diff, G^\diff]$}, and thus the normal subgroup of $G$ generated by $[G^\diff, G^\diff]$ is the only minimal non trivial normal subgroup of $G$. 
This fact is an immediate consequence of Thurston's trick (see the last paragraph of section~\ref{sec.pis}) and Banyaga's theorem. Indeed
let $g$ be a non trivial element in $G$, choose $h \in G$ and a disc $D$ such that $D$, $g(D)$, $h(D)$ are pairwise disjoint. Thurston's trick shows that  the normal subgroup $N(g)$ generated by $g$ in $G$ contains some non trivial commutator of diffeomorphisms, let us denote it by $\Phi$. By Banyaga's theorem any element of $[G^\diff, G^\diff]$ is a product of conjugates (in $G^\diff$) of $\Phi$ and $\Phi^{-1}$, and thus $[G^\diff, G^\diff]$ is included in $N(g)$.

\begin{proof}[Proof of proposition~\ref{prop.profile-diffeo}]
Let $g \in G^\diff$, and fix some smooth isotopy $(g_{t})_{t \in [0,1]}$ from the identity to $g$ in $G^\diff$. Let $M>0$ be such that every trajectory of the isotopy has speed bounded from above by $M$.

We now fix some  $\rho>0$, and let $m$ be the smallest integer such that $m \geq \frac{2}{\rho}$.  We consider some discs $D_{1}, \dots D_{m}$ as in the proof of proposition~\ref{prop.fragmentation-locale}; if we choose the $D_{i}$'s to be  horizontal slices, then for every $i,j$ with $| i -j | > 1$ we get
$$
d(D_{i},D_{j}) > \frac{C}{m}
$$
where $d$ is the euclidean metric of the unit disc and $C$ is some constant (maybe $C=\frac{\pi}{2}$).
Due to the definition of $M$, within any interval of time less than  $\frac{C}{mM}$, no
point moves a distance more than $\frac{C}{m}$:   for every $t,t'\in [0,1]$, for every $x \in \bbD^2$,
 $$
 | t- t'| < \frac{C}{mM} \quad \Longrightarrow \quad    d(g_{t}(x), g_{t'}(x)) < \frac{C}{m}.
 $$
  In particular the topological  disc $g_{t'}g_{t}^{-1}(D_{i})$  remains disjoint from $D_{j}$ for every $|i-j| >1$;
that is,  the diffeomorphism $g_{t'}g_{t}^{-1}$ belongs to the neighbourhood $\cV_{\rho}$ defined in the proof of proposition~\ref{prop.fragmentation-locale}.
Thus we can write $g$ as the product of at most $\frac{mM}{C}+1$ elements of $\cV_{\rho}$,
$$
g =g_{1} =  \left(g_{1} g_{1-\frac{1}{k}}^{-1}\right) \left(g_{1-\frac{1}{k}} g_{1-\frac{2}{k}}^{-1}\right) \cdots \left(g_{\frac{1}{k}} g_{0}^{-1}\right) 
$$
(where $k$ is the integer part of $\frac{mM}{C}+1$).
Each element in $\cV_{\rho}$ is the product of at most $m-1$ elements whose sizes are less than $\rho$, thus we get the estimate
$$
||g||^\diff_{\rho} \leq (m-1) \left(\frac{mM}{C}+1\right).
$$
When $\rho$ tends to $0$, the right-hand side quantity is equivalent to   $\frac{4M}{C}\frac{1}{\rho^2}$, which proves the first point of the proposition.

We turn to the second point.
We first prove the result for some special commutator. Let $D$ be any displaceable disc, say of area $\frac{1}{3}$, 
and let $\Phi$ be any non trivial element of $G^\diff$ supported in the interior of $D$.
Choose some $\Psi\in G^\diff$ such that $\Psi(D)$ is disjoint from $D$.
Let us define $g := [\Phi,\Psi]$. We claim that for any $\rho$ we have  
$$
||g||^\diff_{\rho} \leq \frac{4}{3\rho}.
$$
To prove the claim fix some positive $\rho$.
 It is easy to find, almost explicitly, some $\Psi_{\rho} \in G^\diff$ 
 which is a product of less than $\frac{2}{3\rho}$ elements of size less than $\rho$ and that moves $D$ disjoint from itself (see figure~\ref{fig.moves}).
  Then
 $$
||\ [\Phi,\Psi_{\rho}] \ ||^\diff_{\rho}  \leq ||\Phi\Psi_{\rho} \Phi^{-1}||^\diff_{\rho} + ||\Psi_{\rho}^{-1}||^\diff_{\rho}  = 2||\Psi_{\rho}||^\diff_{\rho} \leq      \frac{4}{3\rho}.
 $$
 We now notice that  the map $g = [\Phi,\Psi]$ is conjugate to  
 $[\Phi,\Psi_{\rho}]$: indeed we may find some $\Theta \in G^\diff$ that is the identity on $D$ and equals $\Psi_{\rho} \Psi^{-1}$ on $\Psi(D)$ (this uses a variation on Lemma~\ref{lem.classical}), and such a $\Theta$ provides the conjugacy.
 Since the fragmentation norm is a conjugacy invariant, this proves the claim.
 \begin{figure}[htbp]
\begin{center}
\def\JPicScale{1.2}
\includegraphics{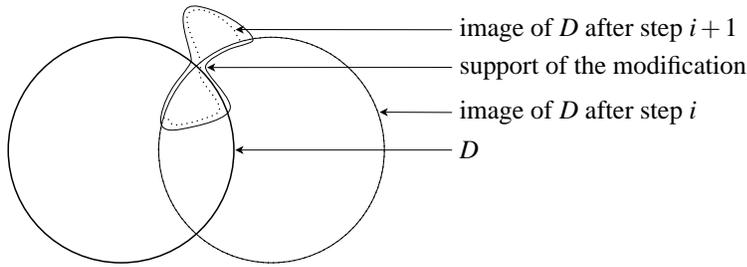}
\caption{How to move a disc disjoint from itself within $\frac{2}{3\rho}$ modifications of size less than $\rho$: after each modification, the area of the intersection of $D$ with its image has decreased by almost~$\frac{\rho}{2}$}
\label{fig.moves}
\end{center}
\end{figure}

We end the proof of the proposition by using Banyaga's theorem. Since the commutator subgroup $[G^\diff, G^\diff]$ is simple, 
 the normal subgroup of $G^\diff$ generated by $g$ in $G^\diff$ is equal to  $[G^\diff, G^\diff]$. As in the $C^0$ case the set of elements $g'\in G^\diff$ satisfying $||g'||^\diff_{\rho} =O(\frac{1}{\rho})$ is a normal subgroup, since it contains $g$ it has to contain $[G^\diff, G^\diff]$. This proves the second point of the proposition.
 \end{proof}

%%%%%%%%%%%%%%%%%%%%%%%%
%%%%%%%%%%%%%%%%%%%%%%%%
\section{Solution of the fragmentation problem for $\rho > \frac{1}{2}$}
\label{sec.quasi-morphisms}
In~\cite{EntovPolterovichPy08} the authors describe some quasi-morphisms on the group of symplectic  diffeomorphisms on various surfaces that are continuous with respect to the $C^0$-topology and, as a consequence, extend continuously to the group of area-preserving homeomorphisms.
 It turns out that their family of quasi-morphisms on the disc provides a solution to the ``easiest'' half of our fragmentation problems:
 properties $(P_{\rho})$ and $(P^\diff_{\rho})$ do not hold when $\rho > \frac{1}{2}$.   We discuss this briefly.

Remember that a map $\phi : G^\diff \to \bbR$ is a \emph{quasi-morphism} if the function
$$
| \phi(gh) - \phi(g) -\phi(h) |
$$
is bounded on $G^\diff \times G^\diff$ by some quantity $\Delta(\phi)$ called the \emph{defect} of $\phi$. A quasi-morphism is called \emph{homogeneous} if it satisfies $\phi(g^n) = n \phi(g)$ for every $g$ and every integer $n$.
The construction of the continuous quasi-morphisms uses the quasi-morphisms of Entov and Polterovich   on the $2$-sphere (\cite{EntovPolterovich03}). This quasi-morphism is a \emph{Calabi} quasi-morphism, that is, it coincides with the Calabi morphism when restricted to those diffeomorphims supported on any displaceable disc. 
By embedding  $\bbD^2$ inside the two-sphere as a non-displaceable\footnote{Note that when the area of the disc tends to the total area of the sphere, the parameter $\rho'$ below tends to $\frac{1}{2}$; the bigger the disc,  the more useful the quasi-morphism, as far as our problem is concerned.}
 disc (\cite{EntovPolterovich03}, Theorem~1.11 and section~5.6), and subtracting the Calabi morphism (\cite{EntovPolterovichPy08}),
one can  get a family $(\phi_{\rho'})_{\rho' \in (\frac{1}{2},1]}$ of homogeneous quasi-morphisms on $G^\diff$. (As these quasi-morphisms extend to $G$, see~\cite{EntovPolterovichPy08}, we could alternatively choose to work with $G$ and $(P_{\rho})$ instead of $G^\diff$ and $(P_{\rho}^\diff)$.) They satisfy the following properties.
\begin{enumerate}
\item $\phi_{\rho'}(g') = 0$ for any $g' \in G^\diff$ whose size is less than $\rho'$,
\item for every $\rho>\rho'$ there exists some $g \in G^\diff$ of size less than $\rho$ with $\phi_{\rho'}(g) \neq 0$.
\end{enumerate}
This entails that for every $\rho \in (\frac{1}{2},1]$, property $(P_{\rho}^\diff)$ do not hold.
  To see this, given some $\rho \in (\frac{1}{2},1]$, we fix some $\rho' \in  (\frac{1}{2},\rho)$,
and we search for some $g \in G^\diff$ having size less than $\rho$ and whose $\rho'$-norm is arbitrarily large (see the easy remarks at the end of section~\ref{sec.fragmentation-metric}).
The first property of $\phi_{\rho'}$
entails, for every $g \in G^\diff$,
$$
\phi_{\rho'}(g) \leq (||g||^\diff_{\rho'}-1)\Delta(\phi_{\rho'}).
$$
On the other hand the second property provides some $g$ with size less than $\rho$ and such that $\phi_{\rho'}(g)\neq0$. Since $\phi_{\rho'}$ is homogeneous, the sequence  $(\phi_{\rho'} (g^n))_{n \geq 0}$  is unbounded.  Thus the sequence  $(||g^n||^\diff_{\rho'})_{n \geq 0}$ is also unbounded, which proves that $(P_{\rho}^\diff)$ do not hold.

\bigskip

We are naturally led to the following question.
\begin{ques}
Does there exist, for every $\rho' \in (0,\frac{1}{2}]$, some homogeneous quasi-morphism $\phi_{\rho'}$ satisfying properties 1 and 2 as expressed above?
\end{ques}
A positive answer would imply a negative answer to Question~\ref{ques.fathi}.

%%%%%%%%%%%%%%%%%%%%%%%%
%%%%%%%%%%%%%%%%%%%%%%%%
\section{Some more remarks}
\label{sec.remarks}
%%%%%%%%%%%%%%%%%%%%%%%

\subsection*{``Lots of'' normal subgroups (if any!)}
The proof of Lemma~\ref{lem.npins}
can  easily be modified to show that, if none of the properties $(P_{\rho})$ holds, then there exists an uncountable family $\cF$ of functions $\varphi$ such that the corresponding family of normal subgroups $(N_{\varphi})_{\varphi \in \cF}$ is totally ordered by inclusion.
The following is another attempt to express that if $G$ is not simple, then it has to contain ``lots of'' normal subgroups.
\begin{coro}
Assume $G$ is not simple. Then every compact subset $K$ of $G$ is included in a proper normal subgroup of $G$.
\end{coro}
Note that the situation is radically different for the diffeomorphisms group $G^\diff$,
since (by Banyaga's theorem~\cite{Banyaga97}, and since the centralizer of $G^\diff$ is trivial) any one-parameter subgroup of diffeomorphisms that is not included in the commutator subgroup $[G^\diff,G^\diff]$
 normally generates $G^\diff$. However  these are not purely algebraic statements since they involve the topology of the groups $G$ and $G^\diff$.
\begin{proof}[Proof of the corollary]
Consider some $\rho \in (0,1]$. Let $\cV_{\rho}$ be the neighbourhood of the identity given by  proposition~\ref{prop.fragmentation-locale}: we have $||g||_{\rho} < \frac{2}{\rho}$ for every $g \in \cV_{\rho}$.
By compactness we may find a finite family $g_{1}, \dots , g_{k}$ such that the sets $g_{i}.\cV_{\rho}$ cover $K$. Thus the fragmentation is also uniform on $K$, in other words the set $K$ is bounded with respect to the norm $||g||_{\rho}$.
Define 
$$
\varphi_{K}(\rho):=\sup\left\{||g||_{\rho}, g \in K \right\}.
$$
This defines a non-increasing function, and clearly $K$ is included in the normal subgroup 
$N_{\varphi_{K}}$. According to Theorem~\ref{theo.good-candidates}, if $G$ is not simple then $N_{\varphi_{K}}$ is a proper subgroup of $G$, which completes the proof of the corollary.
\end{proof}

%%%%%%%%%%%%%%%%%%%%%%%
\subsection*{Other surfaces}
Let $S$ be any compact surface equipped with an area form. Consider the group of homeomorphisms that preserves the measure associated to the area form, and denote by  $G_{0}(S)$  the (normal) subgroup generated by  the homeomorphisms that are supported inside a topological disc. 
For example, $G_{0}(\bbS^2)$ is just the group of orientation and area preserving homeomorphisms of the sphere, and 
 $G_{0}(\bbT^2)$ is the group of orientation and area preserving homeomorphisms of the torus with zero mean rotation vector. The group $G_{0}(S)$ may also be seen as the closure of the group of hamiltonian diffeomorphisms of $S$ inside the group of homeomorphisms.  For every surface $S$, it is an open question whether the  group $G_{0}(S)$ is simple or not.

Exactly as before, on $G_{0}(S)$ we may define the size of an element supported in a topological disc, the family of  fragmention norms $||.||^S_{\rho}$, the family of normal subgroups $N^S_{\varphi}$, and the fragmentation properties $(P^S_{\rho})$.
Then Lemma~\ref{lem.npins} still holds, with the same proof: the failure of all the properties  $(P^S_{\rho})$ would entail that every normal subgroup $N^S_{\varphi}$ is proper, and that $G_{0}(S)$ is not simple. On the other hand the  proof of Lemma~\ref{lem.pis} shows that each property  $(P^{\bbD^2}_{\rho})$ forces the simplicity of $G_{0}(S)$ (alternatively, we could use the original Lemma~\ref{lem.pis} and Thurston trick to see that the simplicity of $G=G_{0}(\bbD^2)$ implies that of $G_{0}(S)$). But property  $(P^{S}_{\rho})$ is \emph{a priori} weaker than  $(P^{\bbD^2}_{\rho})$, because on a general surface one has more space than in the disc to perform the fragmentation.
This prevents us from fully translating the simplicity question into a fragmentation problem on the other surfaces.

Using a different approach, one  might hope to recover this equivalence by adapting to the $C^0$ context the homology machinery introduced by Thurston in the smooth category (see~\cite{Banyaga97} or~\cite{Bounemoura08}, section~2.2).

\end{document}